\documentclass[10pt]{amsart}

\usepackage{amssymb,amsmath,amsthm,amsfonts}

\theoremstyle{plain}
\newtheorem{Theorem}{Theorem}[section]
\newtheorem{Proposition}[Theorem]{Proposition}
\newtheorem{Conjecture}[Theorem]{Conjecture}
\newtheorem{Lemma}[Theorem]{Lemma}

\theoremstyle{definition}

\newtheorem{Example}[Theorem]{Example}
\newtheorem*{Remark}{Remark}

\begin{document}
\title{A finitely presented torsion-free simple group}
\author{Diego Rattaggi}
\address{ETH Z\"urich, Departement Mathematik, R\"amistrasse 101, 8092 Z\"urich, Switzerland}
\email{rattaggi@math.ethz.ch}
\date{April 20, 2004}
\begin{abstract}
We construct a finitely presented torsion-free simple group $\Sigma_0$, acting cocompactly on a product of 
two regular trees. An infinite family of such groups has been introduced by Burger-Mozes (\cite{BM,BMII}).
We refine their methods and get $\Sigma_0$ as an index $4$ subgroup of a group 
$\Sigma < \mathrm{Aut}(\mathcal{T}_{12}) \times \mathrm{Aut}(\mathcal{T}_{8})$ presented by
$10$ generators and $24$ short relations.
For comparison, the smallest virtually simple group of \cite[Theorem~6.4]{BMII}
needs more than $18000$ relations, and the smallest simple group constructed in \cite[Section~6.5]{BMII} needs
even more than $360000$ relations in any finite presentation.
\end{abstract}
\maketitle

\setcounter{section}{-1}
\section{Introduction} \label{SectionIntro}
Burger-Mozes have constructed in \cite{BM,BMII} the first examples of groups which are
simultaneously finitely presented, torsion-free and simple.
Moreover, they are CAT(0), bi-automatic, and have finite cohomological dimension.
These groups can be realized in various ways: as fundamental groups of finite square complexes,
as cocompact lattices in a product of automorphism groups of regular trees 
$\mathrm{Aut}(\mathcal{T}_{2m}) \times \mathrm{Aut}(\mathcal{T}_{2n})$ for sufficiently large $m,n \in \mathbb{N}$,
or as amalgams of finitely generated free groups.
The groups of Burger-Mozes have positively answered several open questions: 
for example Neumann's question (\cite{Neumann}) on the existence of simple amalgams of finitely generated free groups,
or a question of G.~Mess (see \cite[Problem 5.11 (C)]{Kirby}) on the existence 
of finite aspherical complexes with simple fundamental group.
The construction is based on a ``normal subgroup theorem'' (\cite[Theorem~4.1]{BMII}) 
which shows for a certain class of irreducible lattices
acting on a product of trees,
that any non-trivial normal subgroup has finite index.
This statement and its remarkable proof are adapted from the famous analogous theorem of Margulis (\cite[Theorem~IV.4.9]{Margulis}) in the context of
irreducible lattices in higher rank semisimple Lie groups. 
Another important ingredient in the construction of Burger-Mozes is a sufficient criterion (\cite[Proposition~2.1]{BMII})
for the non-residual finiteness of groups acting on a product of trees.  
Even the bare existence of such non-residually finite groups is remarkable, since for example
finitely generated linear groups, or cocompact lattices in $\mathrm{Aut}(\mathcal{T}_{k})$ are always residually finite.
The non-residually finite groups of Burger-Mozes coming from their criterion always do have
non-trivial normal subgroups of infinite index, but appropriate embeddings into groups satisfying the
normal subgroup theorem immediately lead to virtually simple groups.
Unfortunately, these groups and their simple subgroups have very large finite presentations.
We therefore modify the constructions by taking a small non-residually finite group of Wise (\cite[Section~II.5]{Wise}),
embed it into a group $\Sigma < \mathrm{Aut}(\mathcal{T}_{12}) \times \mathrm{Aut}(\mathcal{T}_{8})$
satisfying the normal subgroup theorem, and detect a simple subgroup $\Sigma_0 < \Sigma$ of index~$4$.
Several \textsf{GAP}-programs (\cite{GAP}) have enabled us to find very quickly the groups $\Sigma$ and $\Sigma_0$.
The \textsf{GAP}-code of our programs is documented in 
\cite[Appendix~B]{Rattaggi} for the interested reader. 

\section{Preliminaries} \label{SectionPrel}
As mentioned in the introduction,
the finitely presented torsion-free simple groups of Burger-Mozes and of this paper
appear in various forms.
Probably the most comprehensible approach is to see them as finite index subgroups of
fundamental groups of certain $2$-dimensional cell complexes which are called 
$1$-vertex VH-T-square complexes in \cite{BMII}, complete squared VH-complexes with one vertex in \cite{Wise},
or $(2m,2n)$--complexes in \cite{Rattaggi}. As in \cite{Rattaggi}, we will call these fundamental groups
$(2m,2n)$--groups here.
Let us briefly recall their definition and some properties needed in the construction of the simple example $\Sigma_0$.
Fix $m, n \in \mathbb{N}$ and let $X$ be a finite $2$-dimensional cell complex satisfying the following conditions:
\begin{itemize}
\item Its $1$-skeleton $X^{(1)}$ consists of a single vertex $x$ and oriented loops $a_1^{\pm 1}, \ldots, a_m^{\pm 1}$, 
$b_1^{\pm 1}, \ldots, b_n^{\pm 1}$.
\item There are exactly $mn$ geometric $2$-cells attached to $X^{(1)}$. They are squares with oriented boundary of the form
$aba'b'$, where $a,a' \in A := \{a_1, \ldots, a_m \}^{\pm 1}$ and $b,b' \in B := \{b_1, \ldots, b_n \}^{\pm 1}$.
We think of the elements in $A$ as ``horizontal'' edges and the elements in $B$ as ``vertical'' edges,
and do not distinguish between squares with boundary $aba'b'$, $a'b'ab$, 
$a^{-1} {b'}^{-1} {a'}^{-1} b^{-1}$ and ${a'}^{-1} b^{-1} a^{-1} {b'}^{-1}$, since they induce the same relations in the fundamental group of $X$.
\item The link of the vertex $x$ in $X$ is the complete bipartite graph $K_{2m,2n}$ with $2m + 2n$ vertices
(where the bipartite structure is induced by the decomposition $A \sqcup B$ of $X^{(1)}$ into $2m$ horizontal and $2n$ vertical edges).
In other words, to any pair $(a,b) \in A \times B$ there is a uniquely determined pair $(a',b') \in A \times B$ such that $aba'b'$ is the boundary of
one of the $mn$ squares in $X$.
\end{itemize}
These conditions imply that the universal covering space $\tilde{X}$ of $X$ is a product of two trees $\mathcal{T}_{2m} \times \mathcal{T}_{2n}$,
where $\mathcal{T}_{k}$ denotes the $k$-regular tree.
The fundamental group $\Gamma := \pi_1(X,x)$ of $X$ is called a \emph{$(2m,2n)$--group}.
By construction, it has a finite presentation $\Gamma = \langle a_1, \ldots, a_m, b_1, \ldots, b_n \mid R_{\Gamma} \rangle$,
where $R_{\Gamma}$ consists of $mn$ relations of the form $aba'b' = 1$ induced from the $mn$ squares of $X$,
and $\Gamma$ acts freely and transitively on the vertices of $\mathcal{T}_{2m} \times \mathcal{T}_{2n}$.
Moreover, it follows from the non-positive curvature of $\tilde{X}$ that $\Gamma$ is torsion-free 
(see \cite[Theorem~4.13(2)]{BrHa}).
Equipping $\mathrm{Aut}(\mathcal{T}_{k})$ with the usual topology of simple convergence and 
$\mathrm{Aut}(\mathcal{T}_{2m}) \times \mathrm{Aut}(\mathcal{T}_{2n})$ with the product topology, 
$\Gamma$ can be seen as a cocompact lattice in $\mathrm{Aut}(\mathcal{T}_{2m}) \times \mathrm{Aut}(\mathcal{T}_{2n})$.
We denote by $\mathrm{pr}_1$ and $\mathrm{pr}_2$ the projections of $\Gamma$ to the first and second factor of
$\mathrm{Aut}(\mathcal{T}_{2m}) \times \mathrm{Aut}(\mathcal{T}_{2n})$, respectively, and
let $H_i$, $i = 1,2$, be the closure $\overline{\mathrm{pr}_i(\Gamma)}$.
Fix a vertex $x_h$ of $\mathcal{T}_{2m}$. For each $k \in \mathbb{N}$, we can associate to a $(2m,2n)$--group $\Gamma$
a finite permutation group $P_h^{(k)}(\Gamma) < S_{2m \cdot (2m-1)^{k-1}}$
which describes the action of $\mathrm{Stab}_{H_1}(x_h)$ on the $k$-sphere around $x_h$ in $\mathcal{T}_{2m}$.
These ``local groups'' (or at least their $n$ generators in $S_{2m \cdot (2m-1)^{k-1}}$) 
can be directly computed, given the $mn$ squares of $X$, see \cite[Chapter~1]{BMII} or 
\cite[Section~1.4]{Rattaggi}
for details.
Analogously, one defines local vertical permutation groups $P_v^{(k)}(\Gamma) < S_{2n \cdot (2n-1)^{k-1}}$,
taking the projection to the second factor $\mathrm{Aut}(\mathcal{T}_{2n})$.

There are several equivalent ways to introduce the notion of ``irreducibility'' for
$(2m,2n)$--groups $\Gamma$. For example, $\Gamma$ is called \emph{irreducible} if and only if 
$\mathrm{pr}_2(\Gamma) < \mathrm{Aut}(\mathcal{T}_{2n})$ is not discrete.
Very useful for our purposes is the following criterion of Burger-Mozes, a direct consequence of \cite[Proposition 1.3]{BMII} and
\cite[Proposition 5.2]{BMII}.

\begin{Proposition} \label{PropIrred}
(Burger-Mozes, see also \cite[Proposition~1.2(1b)]{Rattaggi})
Let $\Gamma$ be a $(2m,2n)$--group
such that $n \geq 3$.
Suppose that $P_v^{(1)}(\Gamma)$ is the alternating group $A_{2n}$. Then $\Gamma$ is irreducible if and only if 
$| P_v^{(2)}(\Gamma)| =  | A_{2n} | \cdot | A_{2n-1} |^{2n}$.
\end{Proposition}

Given a $(2m,2n)$--group by its presentation $\Gamma = \langle a_1, \ldots , a_m, b_1, \ldots , b_n \mid R_{\Gamma} \rangle$,
we define a normal subgroup $\Gamma_0$ of index $4$ as kernel of the surjective homomorphism
\begin{align}
\Gamma &\to \mathbb{Z} / 2 \mathbb{Z} \times \mathbb{Z} / 2 \mathbb{Z} \notag \\
a_1, \ldots, a_m    &\mapsto (1 + 2 \mathbb{Z}, 0 + 2 \mathbb{Z}),  \notag \\
b_1, \ldots, b_n    &\mapsto (0 + 2 \mathbb{Z}, 1 + 2 \mathbb{Z}). \notag
\end{align}
Geometrically, $\Gamma_0$ can be seen as fundamental group of a square complex with $4$ vertices,
a $4$-fold regular covering of $X$. The subscript ``$0$'' will always refer to this specific subgroup.

We write $G^{\ast}$ for the intersection of all finite index normal subgroups of a group $G$.
Note that $G^{\ast}$ is a normal subgroup of $G$ and
recall that $G$ is called \emph{residually finite} if and only if $G^{\ast}$ is the trivial group.
It does not matter if one takes the intersection of all finite index subgroups,
or of all finite index \emph{normal} subgroups, as seen in the following elementary lemma.

\begin{Lemma} \label{LemmaNormal}
Let $G$ be a group and $H < G$ a subgroup of finite index $[G : H]$. Then
there is a group $N < H$ such that $N \lhd G$ and $[G : N] \leq [G : H]! < \infty$, in particular
$G^{\ast}$ is also the intersection of all finite index subgroups of $G$.
\end{Lemma}

\begin{proof}
Let $k$ be the finite index $[G : H]$ and write $G$ as a disjoint finite union of left cosets
\[
G = \bigsqcup_{i=1}^{k} g_i H \, .
\]
Left multiplication $g_i H \mapsto g g_i H$ induces a homomorphism $\phi: G \to S_k$ such that
$N := \mathrm{ker} \, \phi < H$ and $[G : N] \leq |S_k| = [G : H]! < \infty$.
\end{proof}

We use the notation $\langle \! \langle g \rangle \! \rangle_{G}$ to denote the \emph{normal closure}
of the element $g \in G$, i.e.\ the intersection of all normal subgroups of $G$ containing $g$. 

\section{The normal subgroup theorem of Burger-Mozes} \label{SectionNST}
Let $T$, $T_1$, $T_2$ be locally finite trees
and let $\Gamma$ be a $(2m,2n)$--group or more generally a subgroup of
$\mathrm{Aut}(T_1) \times \mathrm{Aut}(T_2)$.
For $i=1,2$, let $H_i = \overline{\mathrm{pr}_i(\Gamma)}$ and 
$H_i^{(\infty)}$ be the intersection of all closed finite index subgroups of $H_i$.
A subgroup $H$ of $\mathrm{Aut}(T)$ is called \emph{locally $\infty$-transitive} if $\mathrm{Stab}_H(x)$
acts transitively on the $k$-sphere around $x$ in $T$ for each vertex 
$x$ of $T$ and each $k \in \mathbb{N}$.

The following statement is the general version of the normal subgroup theorem of Burger-Mozes:
\begin{Theorem} \label{NST}
(\cite[Theorem 4.1]{BMII})
Let $\Gamma < \mathrm{Aut}(T_1) \times \mathrm{Aut}(T_2)$ be a cocompact lattice
such that $H_i$ is locally $\infty$-transitive and $H_i^{(\infty)}$ is of finite index in $H_i$, $i=1,2$.
Then, any non-trivial normal subgroup of $\Gamma$ has finite index.
\end{Theorem}

We will use a special version of Theorem~\ref{NST} which directly follows from the discussion in \cite[Chapter~3]{BMI}
and \cite[Chapter~5]{BMII}:
\begin{Theorem} \label{NSTspecial} 
(Burger-Mozes, see also \cite[Proposition 2.1]{Rattaggi})
Let $\Gamma$ be an irreducible $(2m,2n)$--group such that $P_h^{(1)}(\Gamma)$, $P_v^{(1)}(\Gamma)$ are $2$-transitive,
and the stabilizers $\mathrm{Stab}_{P_h^{(1)}(\Gamma)}(\{1\})$, $\mathrm{Stab}_{P_v^{(1)}(\Gamma)}(\{1\})$ 
are non-abelian finite simple groups.
Then any non-trivial normal subgroup of $\Gamma$ has finite index.
\end{Theorem}
We can apply Theorem~\ref{NSTspecial} for example to a group 
$\Lambda < \mathrm{Aut}(\mathcal{T}_6) \times \mathrm{Aut}(\mathcal{T}_6)$,
acting ``locally like $A_6$''.
\begin{Example}
Let
\[
R_{\Lambda} := \left\{ \begin{array}{l l l}
 a_1 b_1 a_1^{-1} b_1^{-1},      &a_1 b_2 a_1^{-1} b_3^{-1},  &a_1 b_3 a_2 b_2^{-1}, \\
                                 &                            &                     \\
 a_1 b_3^{-1} a_3^{-1} b_2,      &a_2 b_1 a_3^{-1} b_2^{-1},  &a_2 b_2 a_3^{-1} b_3^{-1}, \\
                                 &                            &                     \\
 a_2 b_3 a_3^{-1} b_1,           &a_2 b_3^{-1} a_3 b_2,       &a_2 b_1^{-1} a_3^{-1} b_1^{-1} 
\end{array} \right\}
\]  
and $\Lambda := \langle a_1, a_2, a_3, b_1, b_2, b_3 \mid R_{\Lambda} \rangle$ the corresponding $(6,6)$--group.
\end{Example}

\begin{Proposition}
Any non-trivial normal subgroup of $\Lambda$ has finite index.
\end{Proposition}
\begin{proof}
We compute
\[
P_h^{(1)}(\Lambda) = \langle (2,3)(4,5), (1,5,4,2,3), (2,3,5,4,6) \rangle \cong A_6, 
\]
\[
P_v^{(1)}(\Lambda) = \langle (2,3)(4,5), (1,6,3,2)(4,5), (1,4,5,6)(2,3) \rangle \cong A_6, 
\]
and $|P_v^{(2)}(\Lambda)| = 360 \cdot 60^6$.
It follows from Proposition~\ref{PropIrred} that $\Lambda$ is irreducible.
Then we apply Theorem~\ref{NSTspecial}, using that $\mathrm{Stab}_{A_6}(\{1\}) \cong A_5$ is non-abelian simple.
\end{proof}

Computational experiments on finite index subgroups of $\Lambda$ 
(for example using quotpic \cite{quotpic}) lead to the following conjecture:
\begin{Conjecture}
The subgroup $\Lambda_0 < \Lambda$ is simple.
\end{Conjecture}

\section{The simple group $\Sigma_0$} \label{SectionSimple}
The $(8,6)$--group $\Delta$ of Example~\ref{ExD} has been constructed by Wise (\cite{Wise})
to give the first examples of non-residually finite groups in the following three important
classes: finitely presented small cancellation groups, automatic groups,
and groups acting properly discontinuously and cocompactly on CAT(0)-spaces.
We embed $\Delta$ in a $(12,8)$--group $\Sigma$ such that $\Sigma$ has no non-trivial normal subgroups of
infinite index. The explicit knowledge of an element in $\Delta^{\ast}$ enables us to prove that the subgroup
$\Sigma_0 < \Sigma$ is simple.
\begin{Example} \label{ExD}
(See \cite[Section~II.5]{Wise} where $\Delta$ is called $D$)
Let
\[
R_{\Delta} := \left\{ \begin{array}{l l l l}
 a_1 b_1 a_2^{-1} b_2^{-1},
&a_1 b_2 a_1^{-1} b_1^{-1},
&a_1 b_3 a_2^{-1} b_3^{-1},
&a_1 b_3^{-1} a_2^{-1} b_2, \\
&&&\\
a_1 b_1^{-1} a_2^{-1} b_3,
&a_2 b_2 a_2^{-1} b_1^{-1}, 
&a_3 b_1 a_4^{-1} b_2^{-1},
&a_3 b_2 a_3^{-1} b_1^{-1}, \\
&&&\\
a_3 b_3 a_4^{-1} b_3^{-1},
&a_3 b_3^{-1} a_4^{-1} b_2,
&a_3 b_1^{-1} a_4^{-1} b_3,
&a_4 b_2 a_4^{-1} b_1^{-1}
\end{array} \right\}
\]
and $\Delta := \langle a_1, a_2, a_3, a_4, b_1, b_2, b_3 \mid R_{\Delta} \rangle$ the corresponding $(8,6)$--group.
\end{Example}

\begin{Proposition} \label{PropNRFelement}
(\cite[Main Theorem~II.5.5]{Wise})
The group $\Delta$ is non-residually finite and $a_2 a_1^{-1} a_3 a_4^{-1} \in \Delta^{\ast}$.
\end{Proposition}
Observe that $\Delta$ has non-trivial normal subgroups of infinite index,
for example the commutator subgroup $[\Delta, \Delta]$ with infinite quotient 
$\Delta / [\Delta, \Delta] \cong \mathbb{Z} \times \mathbb{Z} \times \mathbb{Z}$.
Our strategy is to embed $\Delta$ as a subgroup in a $(2m,2n)$--group which satisfies the assumptions of Theorem~\ref{NSTspecial},
and to apply the following basic lemma.
\begin{Lemma} \label{LemmaNRF}
Let $G$ be a group and $H < G$ a subgroup. Then $H^{\ast} < G^{\ast}$.
In particular, if $H$ is non-residually finite, then also $G$ is non-residually finite.
\end{Lemma}
\begin{proof}
Let $h \in H^{\ast}$ and $N \lhd G$ any normal subgroup of finite index.
It follows that $N \cap H \lhd G \cap H = H$
such that the index $[H : (N \cap H)] \leq [G : N]$ is finite.
Therefore, $h \in N \cap H < N$. 
\end{proof}

\begin{Example}
Let
\[
R_{\Sigma} := R_{\Delta} \cup
\left\{ \begin{array}{l l l l}
    a_1 b_4 a_3 b_4,
   &a_1 b_4^{-1} a_2 b_4^{-1}, 
   &a_2 b_4 a_5 b_4,  
   &a_3 b_4^{-1} a_4^{-1} b_4^{-1}, \\
&&&\\ 
    a_4 b_4^{-1} a_5 b_4^{-1}, 
   &a_5 b_1 a_6^{-1} b_2,
   &a_5 b_2 a_6^{-1} b_2^{-1}, 
   &a_5 b_3 a_5^{-1} b_3^{-1}, \\
&&&\\ 
    a_5 b_2^{-1} a_6^{-1} b_1^{-1}, 
   &a_5 b_1^{-1} a_6^{-1} b_1, 
   &a_6 b_3 a_6^{-1} b_4^{-1}, 
   &a_6 b_4 a_6^{-1} b_3
\end{array} \right\}
\]
and $\Sigma := \langle a_1, a_2, a_3, a_4, a_5, a_6, b_1, b_2, b_3, b_4 \mid R_{\Sigma} \rangle$ the corresponding $(12,8)$--group.
\end{Example}

\begin{Theorem} \label{ThSimple}
The group $\Sigma_0$ is finitely presented, torsion-free and simple.
\end{Theorem}

\begin{proof}
Being a finite index subgroup of the $(12,8)$--group $\Sigma$, it is clear that
$\Sigma_0$ is finitely presented and torsion-free. It remains to prove that $\Sigma_0$ is simple.

First we show that $\Sigma_0$ has no proper subgroups of finite index.
By construction, $R_{\Sigma}$ contains all twelve elements of
$R_{\Delta}$, hence by \cite[Proposition~II.4.14(1)]{BrHa}, 
this embedding induces an injection on the level of fundamental groups,
i.e.\ $\Delta$ is a subgroup of $\Sigma$.
Let $w := a_2 a_1^{-1} a_3 a_4^{-1} \in \Delta < \Sigma$. By Proposition~\ref{PropNRFelement} and Lemma~\ref{LemmaNRF},
$\Sigma$ is non-residually finite such that $w \in \Sigma^{\ast} \lhd \Sigma$, hence 
$\langle \! \langle w \rangle \! \rangle_{\Sigma} < \Sigma^{\ast}$ by definition of the normal closure.
By a coset enumeration, a computer algebra system like \textsf{GAP} (\cite{GAP}) immediately shows that
adding the relation $w = 1$ to the presentation of $\Sigma$ leads to a finite group of order $4$ 
(the group $\mathbb{Z} / 2\mathbb{Z} \times \mathbb{Z} / 2\mathbb{Z}$), 
in other words $[\Sigma : \langle \! \langle w \rangle \! \rangle_{\Sigma} ] = 4$.
It follows by definition of $\Sigma^{\ast}$ that $\Sigma^{\ast} < \langle \! \langle w \rangle \! \rangle_{\Sigma}$,
thus we have $\Sigma^{\ast} = \langle \! \langle w \rangle \! \rangle_{\Sigma}$.
Since $\Sigma_0$ is a normal subgroup of $\Sigma$ of index $4$, and $w \in \Sigma_0$,
we also get $\langle \! \langle w \rangle \! \rangle_{\Sigma} = \Sigma_0$.
Now it is easy to see that the group $\Sigma_0 = \langle \! \langle w \rangle \! \rangle_{\Sigma} = \Sigma^{\ast}$
has no proper subgroups of finite index as follows:
Assume that $H$ is a finite index subgroup of $\Sigma^{\ast}$, then $H$ has finite index in $\Sigma$
and by Lemma~\ref{LemmaNormal} there is a finite index normal subgroup $N$ of $\Sigma$ such that $N < H < \Sigma^{\ast}$.
By definition of $\Sigma^{\ast}$ we have $\Sigma^{\ast} < N$, hence $N = H = \Sigma^{\ast} = \Sigma_0$.

Next we show that $\Sigma_0$ has no non-trivial normal subgroups of infinite index.
First, we observe that $\Sigma$ is irreducible. 
This is a direct consequence of the fact that $\Sigma$ is non-residually finite,
since reducible $(2m,2n)$--groups are virtually a direct product of two free groups. 
Alternatively, we compute that $P_v^{(2)}(\Sigma)$ has order $20160 \cdot 2520^8$ 
and apply Proposition~\ref{PropIrred}, using 
\begin{align}
P_v^{(1)}(\Sigma) = \langle &(1,2)(4,5)(6,8,7), (1,2,3)(4,5)(7,8), (1,2)(4,5)(6,8,7), \notag \\ 
&(1,2,3)(4,5)(7,8), (1,7)(4,5), (2,8)(3,5,6,4) \rangle \cong A_8. \notag
\end{align}
We also compute
\begin{align}
P_h^{(1)}(\Sigma) = \langle &(5,6)(7,8)(9,10)(11,12), (1,2)(3,4)(5,6)(7,8), \notag \\
&(1,2)(3,4)(9,10)(11,12), (1,11,5,9,10)(2,12,3,4,8) \rangle \cong M_{12}, \notag 
\end{align}
the Mathieu group which acts $5$-transitive on the set $\{ 1, \ldots, 12 \}$. 
Its stabilizer $\mathrm{Stab}_{M_{12}}(\{1\})$ is isomorphic to the non-abelian simple group $M_{11}$ of order $7920$.
By Theorem~\ref{NSTspecial}, any non-trivial normal subgroup of $\Sigma$ has finite index.
Moreover, applying Theorem~\ref{NST}, any non-trivial normal subgroup of $\Sigma_0 = \Sigma^{\ast}$ has finite index.
Note that one uses here again the fact that $\Sigma^{\ast}$ has finite index in $\Sigma$ 
(see the reasoning leading to \cite[Corollary 5.4]{BMII}).
\end{proof}

\begin{Remark}
The simple group $\Sigma_0$ has amalgam decompositions of the form $F_7 \ast _{F_{73}} F_7$ and
$F_{11} \ast _{F_{81}} F_{11}$, where $F_k$ denotes the free group of rank $k$. 
This follows from \cite[Theorem~I.1.18]{Wise},
see also \cite[Proposition~1.4]{Rattaggi}.
The smallest candidate for being a finitely presented torsion-free simple group
in the construction of virtually simple groups in \cite[Theorem~6.4]{BMII} has amalgam decompositions
$F_{349} \ast _{F_{75865}} F_{349}$ and $F_{217} \ast _{F_{75601}} F_{217}$.
The amalgam decompositions of the smallest simple group constructed in \cite[Theorem~6.5]{BMII} 
are $F_{7919} \ast _{F_{380065}} F_{7919}$ and $F_{47} \ast _{F_{364321}} F_{47}$.
\end{Remark}

\begin{Remark}
It is easy to get an explicit finite presentation of $\Sigma_0$: 
Either we can take the fundamental group of the covering space corresponding to the subgroup $\Sigma_0 < \Sigma$,
or we take a presentation of an amalgam mentioned in the remark above 
(note that its explicit construction also makes use of this covering space and additionally the Seifert-van Kampen Theorem).
A third possibility is to use a computer algebra system like \textsf{GAP} (\cite{GAP}),
which has implemented a Reidemeister-Schreier method. Applying this method and Tietze transformations to
reduce the number of generators, 
we get a presentation of $\Sigma_0$ with $3$ generators and $62$ relations of total length $4866$.
\end{Remark}

\section{Generalization}
The proof of Theorem~\ref{ThSimple} shows that if we embed the non-residually finite $(8,6)$--group $\Delta$ 
in a $(2m,2n)$--group $\Gamma$ such that
$P_h^{(1)}(\Gamma)$, $P_v^{(1)}(\Gamma)$ are $2$-transitive
and $\mathrm{Stab}_{P_h^{(1)}(\Gamma)}(\{1\})$, 
$\mathrm{Stab}_{P_v^{(1)}(\Gamma)}(\{1\})$ are non-abelian simple,
then the normal subgroup $\langle \! \langle a_2 a_1^{-1} a_3 a_4^{-1} \rangle \! \rangle_{\Gamma}$ has finite index in $\Gamma$,
and $\Gamma^{\ast} = \langle \! \langle a_2 a_1^{-1} a_3 a_4^{-1} \rangle \! \rangle_{\Gamma}$ is
a finitely presented torsion-free simple group.
In this way, we have constructed many more such simple groups $\Gamma^{\ast}$
for 
\[
(2m,2n) \in \{ (10,10), (10,12), (12,8), (12,10), (12,12) \}, 
\]
see Table~\ref{TableSimple}.
In this table, $D_k$ denotes the dihedral group of order $2k$.
Note that the index $[\Gamma : \Gamma^{\ast}]$ can be larger than $4$,
and that we have $[\Gamma, \Gamma] = \Gamma^{\ast}$ in most cases of Table~\ref{TableSimple}.
\newpage
\begin{table}[ht]
\begin{center}
\begin{tabular}{| l | l | l | c |}
\hline
$P_h^{(1)}(\Gamma)$ & $P_v^{(1)}(\Gamma)$ & $\Gamma / \Gamma^{\ast}$ & $\Gamma / [\Gamma, \Gamma ]$ \\ \hline \hline
$A_{10}$  & $A_{10}$  & $\mathbb{Z} / 2 \mathbb{Z} \times \mathbb{Z} / 2 \mathbb{Z}$ & $=$ \\ \hline
$A_{10}$  & $A_{10}$  & $\mathbb{Z} / 2 \mathbb{Z} \times \mathbb{Z} / 2 \mathbb{Z} \times \mathbb{Z} / 2 \mathbb{Z}$ & $=$ \\ \hline
$A_{10}$  & $A_{10}$  & $\mathbb{Z} / 2 \mathbb{Z} \times \mathbb{Z} / 4 \mathbb{Z}$ & $=$ \\ \hline
$A_{10}$  & $A_{10}$  & $\mathbb{Z} / 2 \mathbb{Z} \times \mathbb{Z} / 6 \mathbb{Z}$ & $=$ \\ \hline
$A_{10}$  & $A_{10}$  & $\mathbb{Z} / 2 \mathbb{Z} \times \mathbb{Z} / 2 \mathbb{Z} \times \mathbb{Z} / 4 \mathbb{Z}$ & $=$ \\ \hline
$A_{10}$  & $A_{10}$  & $\mathbb{Z} / 2 \mathbb{Z} \times \mathbb{Z} / 8 \mathbb{Z}$ & $=$ \\ \hline
$A_{10}$  & $A_{10}$  & $\mathbb{Z} / 2 \mathbb{Z} \times \mathbb{Z} / 10 \mathbb{Z}$ & $=$ \\ \hline
$A_{10}$  & $A_{10}$  & $\mathbb{Z} / 2 \mathbb{Z} \times \mathbb{Z} / 2 \mathbb{Z} \times \mathbb{Z} / 6 \mathbb{Z}$ & $=$ \\ \hline
$A_{10}$  & $A_{10}$  & $\mathbb{Z} / 2 \mathbb{Z} \times \mathbb{Z} / 12 \mathbb{Z}$ & $=$ \\ \hline
$A_{10}$  & $A_{10}$  & $\mathbb{Z} / 2 \mathbb{Z} \times \mathbb{Z} / 2 \mathbb{Z} \times \mathbb{Z} / 8 \mathbb{Z}$ & $=$ \\ \hline
$A_{10}$  & $A_{10}$  & $\mathbb{Z} / 2 \mathbb{Z} \times \mathbb{Z} / 20 \mathbb{Z}$ & $=$ \\ \hline \hline
$A_{10}$  & $A_{12}$  & $\mathbb{Z} / 2 \mathbb{Z} \times \mathbb{Z} / 2 \mathbb{Z}$ & $=$ \\ \hline
$A_{10}$  & $A_{12}$  & $D_6$ & $\mathbb{Z} / 2 \mathbb{Z} \times \mathbb{Z} / 2 \mathbb{Z}$ \\ \hline
$A_{10}$  & $A_{12}$  & $\mathbb{Z} / 2 \mathbb{Z} \times \mathbb{Z} / 2 \mathbb{Z} \times \mathbb{Z} / 2 \mathbb{Z}$ & $=$ \\ \hline
$A_{10}$  & $A_{12}$  & $S_3 \times \mathbb{Z} / 2 \mathbb{Z} \times \mathbb{Z} / 2 \mathbb{Z}$ & $\mathbb{Z} / 2 \mathbb{Z} \times \mathbb{Z} / 2 \mathbb{Z} \times \mathbb{Z} / 2 \mathbb{Z}$ \\ \hline
$A_{10}$  & $A_{12}$  & $\mathbb{Z} / 2 \mathbb{Z} \times \mathbb{Z} / 4 \mathbb{Z}$ & $=$ \\ \hline \hline
$A_{12}$  & $A_{8}$  & $\mathbb{Z} / 2 \mathbb{Z} \times \mathbb{Z} / 2 \mathbb{Z}$ & $=$ \\ \hline
$A_{12}$  & $A_{8}$  & $\mathbb{Z} / 2 \mathbb{Z} \times \mathbb{Z} / 4 \mathbb{Z}$ & $=$ \\ \hline
$M_{12}$  & $A_{8}$  & $\mathbb{Z} / 2 \mathbb{Z} \times \mathbb{Z} / 2 \mathbb{Z}$ & $=$ \\ \hline \hline
$A_{12}$  & $A_{10}$  & $\mathbb{Z} / 2 \mathbb{Z} \times \mathbb{Z} / 2 \mathbb{Z}$ & $=$ \\ \hline
$A_{12}$  & $A_{10}$  & $D_6$ & $\mathbb{Z} / 2 \mathbb{Z} \times \mathbb{Z} / 2 \mathbb{Z}$ \\ \hline
$A_{12}$  & $A_{10}$  & $D_5 \times \mathbb{Z} / 2 \mathbb{Z}$ & $\mathbb{Z} / 2 \mathbb{Z} \times \mathbb{Z} / 2 \mathbb{Z}$ \\ \hline
$A_{12}$  & $A_{10}$  & $\mathbb{Z} / 2 \mathbb{Z} \times \mathbb{Z} / 2 \mathbb{Z} \times \mathbb{Z} / 2 \mathbb{Z}$ & $=$ \\ \hline
$A_{12}$  & $A_{10}$  & $\mathbb{Z} / 2 \mathbb{Z} \times \mathbb{Z} / 4 \mathbb{Z}$ & $=$ \\ \hline
$A_{12}$  & $A_{10}$  & $D_4 \times \mathbb{Z} / 2 \mathbb{Z}$ & $\mathbb{Z} / 2 \mathbb{Z} \times \mathbb{Z} / 2 \mathbb{Z} \times \mathbb{Z} / 2 \mathbb{Z}$ \\ \hline
$A_{12}$  & $A_{10}$  & $\mathbb{Z} / 2 \mathbb{Z} \times \mathbb{Z} / 6 \mathbb{Z}$ & $=$ \\ \hline
$A_{12}$  & $A_{10}$  & $\mathbb{Z} / 2 \mathbb{Z} \times \mathbb{Z} / 8 \mathbb{Z}$ & $=$ \\ \hline
$A_{12}$  & $A_{10}$  & $\mathbb{Z} / 2 \mathbb{Z} \times \mathbb{Z} / 10 \mathbb{Z}$ & $=$ \\ \hline
$A_{12}$  & $A_{10}$  & $\mathbb{Z} / 2 \mathbb{Z} \times \mathbb{Z} / 2 \mathbb{Z} \times \mathbb{Z} / 6 \mathbb{Z}$ & $=$ \\ \hline
$M_{12}$  & $A_{10}$  & $\mathbb{Z} / 2 \mathbb{Z} \times \mathbb{Z} / 2 \mathbb{Z}$ & $=$ \\ \hline \hline
$A_{12}$  & $A_{12}$  & $\mathbb{Z} / 2 \mathbb{Z} \times \mathbb{Z} / 2 \mathbb{Z}$ & $=$ \\ \hline
$A_{12}$  & $A_{12}$  & $\mathbb{Z} / 2 \mathbb{Z} \times \mathbb{Z} / 2 \mathbb{Z} \times \mathbb{Z} / 2 \mathbb{Z}$ & $=$ \\ \hline
$A_{12}$  & $A_{12}$  & $\mathbb{Z} / 2 \mathbb{Z} \times \mathbb{Z} / 6 \mathbb{Z}$ & $=$ \\ \hline
\end{tabular}
\end{center}
\caption{List of simple groups $\Gamma^{\ast}$} \label{TableSimple}
\end{table}

\end{document}